\documentclass[12pt,reqno]{amsart}
\usepackage{enumerate, latexsym, amsmath, amsfonts, amssymb, amsthm, color}
\def\pmod #1{\ ({\rm{mod}}\ #1)}
\def\Z{\Bbb Z}
\def\N{\Bbb N}

\def\M{\Bbb M}
\def\l{\left}
\def\r{\right}
\def\bg{\bigg}
\def\({\bg(}
\def\){\bg)}
\def\t{\text}
\def\f{\frac}

\def\ls{\leqslant}

\def\bi{\binom}
\def\al{\alpha}

\def\eq{\equiv}

\def\Proof{\noindent{\it Proof}}
\def\Ack{\medskip\noindent {\bf Acknowledgment}}
\def\M#1#2{{#1\brack #2}}
\theoremstyle{plain}
\newtheorem{theorem}{Theorem}

\newtheorem{lemma}{Lemma}
\newtheorem{corollary}{Corollary}

\theoremstyle{definition}

\theoremstyle{remark}

\makeatletter
\@namedef{subjclassname@2020}{%
  \textup{2020} Mathematics Subject Classification}
\makeatother
 \vspace{4mm}

\begin{document}

\hbox{Bull. Malays. Math. Sci. Soc. 46 (2023), no.\,4, Article No. 119.}
\medskip

\title
[Legendre symbols related to certain determinants]
{Legendre symbols related to \\ certain determinants}

\author
[Xin-Qi Luo and Zhi-Wei Sun] {Xin-Qi Luo and Zhi-Wei Sun}

\address {(Xin-Qi Luo) Department of Mathematics, Nanjing
University, Nanjing 210093, People's Republic of China}
\email{lxq15995710087@163.com}

\address{(Zhi-Wei Sun, corresponding author) Department of Mathematics, Nanjing
University, Nanjing 210093, People's Republic of China}
\email{zwsun@nju.edu.cn}

\keywords{Legendre symbols, determinants, generalized trinomial coefficients, Lucas sequences.
\newline \indent 2020 {\it Mathematics Subject Classification}. Primary 11C20, 15A15; Secondary 11A15, 11B39.}

\begin{abstract}
Let $p$ be an odd prime. For $b,c\in\mathbb Z$, Sun introduced the determinant
$$D_p(b,c)=\left|(i^2+bij+cj^2)^{p-2}\right|_{1\leqslant i,j \leqslant p-1},$$
and investigated the Legendre symbol $(\frac{D_p(b,c)}p)$.
Recently Wu, She and Ni proved that $(\frac{D_p(1,1)}p)=(\frac {-2}p)$ if $p\equiv2\pmod 3$, which confirms a previous conjecture of Sun. In this paper we determine $(\frac{D_p(1,1)}p)$ in the case $p\equiv1\pmod3$.
Sun proved that $D_p(2,2)\equiv0\pmod p$ if $p\equiv3\pmod4$, in contrast we prove that
$(\frac{D_p(2,2)}p)=1$ if $p\equiv1\pmod8$, and $(\frac{D_p(2,2)}p)=0$ if $p\equiv5\pmod8$.
Our tools include generalized trinomial coefficients and Lucas sequences.
\end{abstract}
\maketitle

\section{Introduction}
\setcounter{lemma}{0}
\setcounter{theorem}{0}
\setcounter{corollary}{0}
\setcounter{remark}{0}
\setcounter{equation}{0}

For an $n\times n$ matrix $[a_{ij}]_{1\leqslant i,j\leqslant n}$ over a commutative ring, we use $|a_{ij}|_{1\ls i,j\ls n}$ to denote its determinant.

Let $p$ be an odd prime, and let $(\frac{.}{p})$ be the Legendre symbol.
Carlitz \cite{C} determined the characteristic polynomial of the matrix
$$\left[x+\left(\frac{i-j}{p}\right)\right]_{1\leqslant i,j\leqslant p-1},
$$
and Chapman \cite{C1} evaluated the determinant
$$
\left| x+\left(\frac{i+j-1}{p}\right)\right|_{1\leqslant i,j\leqslant (p-1)/{2}}.
$$
Vsemirnov \cite{V12,V13} confirmed a challenging conjecture of Chapman
by evaluating the determinant
$$\left|\left(\frac{j-i}{p}\right)\right|_{1\leqslant i,j\leqslant (p+1)/{2}}.
$$
Sun \cite{S19} studied some determinants whose entries have the form $(\f{i^2+cij+dj^2}p)$,
where $c,d\in\Z$; in particular he introduced
$$(c,d)_p:=\left|\left(\frac{i^2+cij+dj^2}{p}\right)\right|_{1\ls i,j\ls p-1}$$
and $$[c,d]_p:=\left|\left(\frac{i^2+cij+dj^2}{p}\right)\right|_{0\ls i,j\ls p-1},$$
and proved that if $(\f dp)=1$ then
$$[c,d]_p=\begin{cases}\f{p-1}2(c,d)_p&\t{if}\ p\nmid c^2-4d,
\\\f{1-p}{p-2}(c,d)_p&\t{if}\ p\mid c^2-4d.\end{cases}$$

For any prime $p\eq3\pmod4$, Sun \cite[Remark 1.3]{S19} showed that 
$$\left|\frac{1}{i^2+j^2}\right|_{1\leqslant i,j \leqslant (p-1)/2}\eq\l(\f 2p\r)\pmod p.
$$
For each prime $p\equiv 5\pmod{6}$, Sun \cite{S19} conjectured that $$2\left|\frac{1}{i^2-ij+j^2}\right|_{1\leqslant i,j \leqslant p-1}$$
 is a quadratic residue modulo $p$. This was recently confirmed by  Wu, She, and Ni \cite{WSN}.

Let $p$ be an odd prime. For $b,c\in\Z$,
Sun \cite{S22} investigated the determinant
\begin{equation}D_p(b,c)=\left|(i^2+bij+cj^2)^{p-2}\right|_{1\leqslant i,j \leqslant p-1},
\end{equation}
and studied the Legendre symbol $(\f{D_p(b,c)}p)$. By Fermat's little theorem,
$$(i^2+bij+cj^2)^{p-2}\eq\begin{cases}\f1{i^2+bij+cj^2}\pmod p&\t{if}\ i^2+bij+cj^2\not\eq0\pmod p,
\\0\pmod p&\t{if}\ i^2+bij+cj^2\eq0\pmod p.
\end{cases}$$
As pointed out in \cite[(1.7)]{S22},
$$D_p(-b,c)\eq\l(\f{-1}p\r)D_p(b,c)\pmod p.$$
Thus, in view of the Wu-She-Ni result \cite{WSN}, if $p\eq2\pmod3$ then
$$D_p(1,1)=\l(\f{(-1)^{(p-1)/2}D_p(-1,1)}p\r)=\l(\f{-1}p\r)\l(\f 2p\r)=\l(\f{-2}p\r).$$

Our first purpose is to determine the Legendre symbol $(\f{D_p(1,1)}p)$ for any prime $p\eq1\pmod3$.

\begin{theorem}\label{Th1.1} Let $p$ be a prime with $p\eq1\pmod 3$. Then
\begin{equation}\l(\f{D_p(1,1)}p\r)=\begin{cases}0&\t{if}\ p\eq7\pmod 9,
\\1&\t{otherwise}.\end{cases}
\end{equation}
\end{theorem}

Sun \cite{S22} proved that $D_p(2,2)\eq0\pmod p$ for any prime $p\eq3\pmod4$.
In contrast, we obtain the following result.

\begin{theorem}\label{Th1.2} Let $p$ be a prime with $p\eq1\pmod4$. Then
\begin{equation}\left(\frac{D_p(2,2)}{p}\right)=\begin{cases}1&\t{if}\ p\eq1\pmod8,\\0&\t{if}\ p\eq5\pmod 8.
\end{cases}\end{equation}
\end{theorem}

In the next section, we will provide some lemmas on generalized trinomial coefficients.
We are going to prove Theorems 1.1 and 1.2 in Sections 3 and 4 respectively.

\section{On generalized trinomial coefficients}
\setcounter{lemma}{0}
\setcounter{theorem}{0}
\setcounter{corollary}{0}
\setcounter{remark}{0}
\setcounter{equation}{0}

Let $n\in\N=\{0,1,2,\ldots\}$. The trinomial coefficients $\bi nk_2\ (k=-n,\ldots,n)$ (cf. \cite{A}) are defined by
$$(x+1+x^{-1})^n=\sum_{k=-n}^n\binom{n}{k}_2x^k,$$
and the number $T_n=\bi n0_2$ is called a central trinomial coefficient.

Let $n\in\N$ and $b,c\in\Z$. We define the generalized trinomial coefficients
$$\binom{n}{k}_{b,c}\ \ (k\in\Z)$$ by
\begin{equation}\label{T}
\l(x+b+\frac{c}{x}\r)^n=\sum_{k\in\Z}\binom{n}{k}_{b,c}x^k.
\end{equation}
Obviously $\bi nk_{b,c}=0$ if $|k|>n$.
Note that $\bi n0_{b,c}$ is just the generalized central trinomial coefficient $T_n(b,c)$
studied in \cite{S14a,S14b}. Clearly,
$$\l(x+2+\f 1x\r)^n=\f{(x+1)^{2n}}{x^n}=\sum_{k=-n}^n\bi{2n}{n+k}x^k$$
and thus $\bi nk_{2,1}=\bi{2n}{n+k}$ for all $k=-n,\ldots,n$.
When $c\not=0$, replacing $x$ in \eqref{T} by $c/x$ we get
$$\l(\f cx+b+x\r)^n=\sum_{k\in\Z}\bi nk_{b,c}\l(\f cx\r)^k=\sum_{k\in\Z}\bi n{-k}c^{-k}x^k,$$
and hence
\begin{equation}\label{-k}\bi nk_{b,c}=\bi{n}{-k}_{b,c}c^{-k}\ \ \t{for all}\ k\in\Z
\end{equation}
in view of \eqref{T}.

Let $b,c\in\Z$. For any $n\in\Z^+=\{1,2,3\ldots\}$, as
$(x+b+\f cx)^n$ equals
$$x\l(x+b+\f cx\r)^{n-1}+b\l(x+b+\f cx\r)^{n-1}+\f cx\l(x+b+\f cx\r)^{n-1},$$
we have the recurrence
\begin{equation}\label{recur}\binom{n}{k}_{b,c}=\binom{n-1}{k-1}_{b,c}
+b\binom{n-1}{k}_{b,c}+c\binom{n-1}{k+1}_{b,c}\end{equation}
for any $k\in\Z$.

\begin{lemma}\label{Lem2.1}
Let $p$ be an odd prime, and let $b,c\in\Z$.
For $k\in\{-p+2,\ldots,p-2\}$, we have
\begin{equation}\label{p-2}\begin{aligned}&(4c-b^2)\binom{p-2}{k}_{b,c}
\\\equiv\ &
\begin{cases}\binom{p-1}{-1}_{b,c}+c\binom{p-1}{1}_{b,c}-b\pmod{p}&\t{if}\ k=0,
\\(k+1)\binom{p-1}{k-1}_{b,c}-(k-1)c\binom{p-1}{k+1}_{b,c}\pmod{p}&\t{if}\ 0<|k|\ls p-2.
\end{cases}\end{aligned}\end{equation}
\end{lemma}

\Proof. For the sake of convenience, for $n\in\N$ and $k\in\Z$ we simply write
$\M nk$ for $\bi nk_{b,c}$.

Taking derivatives of both sides of the identity
\begin{equation}\label{k}\sum_{k=-p}^p\M{p}{k}x^k= \l(x+b+\frac{c}{x}\r)^p,
\end{equation}
 we get
\begin{equation}\label{k-1}\sum_{k=-p}^p\M{p}{k}kx^{k-1}
=p\l(x+b+\frac{c}{x}\r)^{p-1}\l(1-\frac{c}{x^2}\r).
\end{equation}
Taking derivatives of both sides of \eqref{k-1}, we obtain
\begin{equation}\begin{aligned}\label{k-2}&\sum_{k=-p}^p\M{p}{k}k(k-1)x^{k-2}
\\=&p(p-1)\l(x+b+\frac{c}{x}\r)^{p-2}\l(1-\frac{c}{x^2}\r)^2+p\l(x+b+\f cx\r)^{p-1}\f{2c}{x^3}.
\end{aligned}\end{equation}
For each $k=-p,\ldots,p$, comparing the coefficients of $x^{k-1}$ on both sides of \eqref{k-1} we get
\begin{equation}\label{(k)}\frac{k}{p}\M{p}{k}=\M{p-1}{k-1}-c\M{p-1}{k+1};
\end{equation}
similarly, comparing the coefficients of $x^{k-2}$ on both sides of \eqref{k-2} we obtain
\begin{equation} \label{k(k-1)}\begin{aligned}&\frac{k(k-1)}{p}\M{p}{k}-2c\M{p-1}{k+1}
\\=\ &(p-1)\left( \M{p-2}{k-2}-2c\M{p-2}{k}+c^2\M{p-2}{k+2}\right).
\end{aligned}\end{equation}

Let $k\in\{-p,\ldots,p\}$.
With the aid of the recurrence \eqref{recur}, we have
\begin{align*}&\M{p-2}{k-2}-2c\M{p-2}k+c^2\M{p-2}{k+2}
\\=\ &\(\M{p-1}{k-1}-b\M{p-2}{k-1}-c\M{p-2}k\)-2c\M{p-2}k
\\\ &+c\(\M{p-1}{k+1}-\M{p-2}k-b\M{p-2}{k+1}\)
\\=\ &\M{p-1}{k-1}+c\M{p-1}{k+1}-4c\M{p-2}k-b\l(\M{p-2}{k-1}+c\M{p-2}{k+1}\r)
\\=\ &\M{p-1}{k-1}+c\M{p-1}{k+1}-4c\M{p-2}k-b\l(\M{p-1}k-b\M{p-2}k\r)
\\=\ &\M{p-1}{k-1}-b\M{p-1}k+c\M{p-1}{k+1}+(b^2-4c)\M{p-2}k
\\=\ &\M{p-1}{k-1}-\(\M pk-\M{p-1}{k-1}-c\M{p-1}{k+1}\)
\\\quad&+c\M{p-1}{k+1}+(b^2-4c)\M{p-2}k
\\=\ &-\M pk+2\M{p-1}{k-1}+2c\M{p-1}{k+1}+(b^2-4c)\M{p-2}k.
\end{align*}
Combining this with \eqref{(k)} and \eqref{k(k-1)}, we get
\begin{align*}
&(k-1)\l(\M{p-1}{k-1}-c\M{p-1}{k+1}\r)-2cp\M{p-1}{k+1}
\\=\ &(1-p)\l(\M pk-2\M{p-1}{k-1}-(b^2-4c)\M{p-2}k\r)
\end{align*}
and hence
\begin{equation}\label{mod}(k+1)\M{p-1}{k-1}-(k-1)c\M{p-1}{k+1}
\eq\M pk-(b^2-4c)\M{p-2}k\pmod p.
\end{equation}

Since
\begin{align*}\sum_{k=-p}^p\M{p}{k}x^{p+k}&=\l(x^2+bx+c\r)^p
\\&\equiv x^{2p}+b^px^p+c^p\equiv x^{2p}+bx^p+c \pmod {p},\end{align*}
we see that
\begin{equation}\label{pk}\M{p}{k}\equiv \begin{cases}b\pmod p&\t{if}\ k=0,
\\1\pmod p&\t{if}\ k=p,
\\c\pmod p&\t{if}\ k=-p,
\\0\pmod p&\t{if}\ k\in\{\pm1,\ldots,\pm(p-1)\}.
\end{cases}
\end{equation}
Combining this with \eqref{mod}, we immediately obtain the desired \eqref{p-2}.
\qed

\begin{lemma}\label{Lem2.2}
Let $p$ be an odd prime, and let $b,c\in\Z$. Then
\begin{equation}\label{(p-2)}
\begin{aligned}&(x^2+bx+c)^{p-2}-c^{p-2}
\\\eq\ & \bi{p-2}1_{b,c}x^{p-1}+\bi{p-2}0_{b,c}x^{p-2}
\\&+\sum_{1<k<p-1}\l(\bi{p-2}k_{b,c}+c^{p-1-k}\bi{p-2}{p-1-k}_{b,c}\r)x^{k-1}\pmod p.
\end{aligned}
\end{equation}
\end{lemma}
\noindent{\bf Proof}. In view of \eqref{-k}, we have
\begin{align*}&(x^2+bx+c)^{p-2}-\bi{p-2}0_{b,c}x^{p-2}
\\=\ &\sum_{k=-(p-2)\atop k\not=0}^{p-2}\bi{p-2}k_{b,c}x^{p-2+k}
\\=\ &\sum_{k=1}^{p-2}\l(\bi{p-2}k_{b,c}x^{p-2+k}+\bi{p-2}{-k}_{b,c}x^{p-2-k}\r)
\\=\ &\sum_{k=1}^{p-2}\l(\bi{p-2}k_{b,c}x^{p-2+k}+\bi{p-2}{k}_{b,c}c^kx^{p-2-k}\r)
\\=\ &\sum_{k=1}^{p-2}\l(\bi{p-2}k_{b,c}x^{p-2+k}+\bi{p-2}{p-1-k}_{b,c}c^{p-1-k}x^{k-1}\r)
\\=\ &\sum_{k=1}^{p-2}\l(\bi{p-2}k_{b,c}x^{p-1}+\bi{p-2}{p-1-k}_{b,c}c^{p-1-k}\r)x^{k-1}.
\end{align*}
Note that
$$\bi{p-2}1_{b,c}x^{p-1}+\bi{p-2}{p-1-1}_{b,c}c^{p-1-1}=\bi{p-2}1_{b,c}x^{p-1}+c^{p-2}.$$
Therefore, from the above we get the desired \eqref{(p-2)}. \qed

For convenience, for an assertion $A$ we set
$$[A]=\begin{cases}1&\t{if}\ A \ \t{holds},
\\0&\t{otherwise}.\end{cases}$$

\begin{corollary}\label{Cor2.1} Let $p$ be an odd prime. If $p\eq1\pmod3$, then
\begin{equation}\label{x^2}\begin{aligned}(x^2+x+1)^{p-2}\eq&\ 1+\f 23x^{p-1}-\f13x^{p-2}
\\&\ +\sum_{k=2}^{p-2}\l(k\l(\f k3\r)+[3\mid k-1]-\f13\r)x^{k-1}\pmod p.
\end{aligned}\end{equation}
\end{corollary}
\Proof.
By \eqref{pk}, we have
$$\bi{p}0_2=1,\ \t{and}\ \bi pk_2=0\ \t{for}\ k=1,\ldots,p-1.$$
Combining this with \eqref{recur}, we see that
$$\bi{p-1}{k-1}_2+\bi{p-1}k_2+\bi{p-1}{k+1}_2\eq0\pmod p$$
for all $k=1,\ldots,p-1.$
In view of this and the easy equalities
$$\bi{p-1}p_2=0\ \ \t{and}\ \ \bi{p-1}{p-1}_2=1,$$
by induction we obtain that
\begin{equation}\label{k/3}\bi{p-1}{p-k}_2\eq\l(\f k3\r)\pmod p\ \ \ \t{for all}\ k=0,1,\ldots,p.
\end{equation}
Note that
\begin{equation}\label{Tp}T_p=\bi{p-1}0_2\eq\l(\f p3\r)3^{p-1}\pmod{p^2}
\end{equation}
as proved by Cao and Sun \cite{CS}.

By Lemma 2.1, \eqref{-k} and \eqref{k/3}, we have
\begin{align*}3\bi{p-2}0_2&\eq\bi{p-1}{-1}_2+\bi{p-1}1_2-1=2\bi{p-1}1_2-1
\\&\eq2\l(\f{p-1}3\r)-1\pmod p.
\end{align*}
Combining Lemma 2.1 and \eqref{k/3}, we see that for each $k=1,\ldots,p-2$ we have
\begin{align*}3\bi{p-2}k_2&\eq(k+1)\bi{p-1}{k-1}_2-(k-1)\bi{p-1}{k+1}_2
\\&\eq(k+1)\l(\f{p-k+1}3\r)-(k-1)\l(\f{p-k-1}3\r)\pmod p
\end{align*}
and
\begin{align*}&3\l(\bi{p-2}k_2+\bi{p-2}{p-1-k}_2\r)
\\\eq\ &(k+1)\l(\f{p-k+1}3\r)-(k-1)\l(\f{p-k-1}3\r)
\\&+(p-1-k+1)\l(\f{p-(p-1-k)+1}3\r)
\\&-(p-1-k-1)\l(\f{p-(p-1-k)-1}3\r)
\\\eq\ &(k+1)\l(\f{p-k+1}3\r)-(k-1)\l(\f{p-k-1}3\r)
\\&-k\l(\f{k+2}3\r)+(k+2)\l(\f k3\r)
\pmod p.
\end{align*}

Now we suppose that $p\eq1\pmod3$. By the last paragraph,
\begin{equation}\label{p-20}\bi{p-2}0_2\eq-\f13\pmod p,
\end{equation}
and for each $k=1,\ldots,p-2$ we have
\begin{align*}&3\l(\bi{p-2}k_2+\bi{p-2}{p-1-k}_2\r)
\\\eq\ & (k+1)\l(\f{-k-1}3\r)-(k-1)\l(\f{-k}3\r)
-k\l(\f{k-1}3\r)+(k+2)\l(\f k3\r)
\\=\ &(2k+1)\l(\f k3\r)-\l(\f{k+1}3\r)
-k\l(\l(\f{k+1}3\r)+\l(\f{k-1}3\r)\r)
\\=\ &(3k+1)\l(\f k3\r)-\l(\f{k+1}3\r)\pmod p.
\end{align*}

Applying Lemma 2.2 with $b=c=1$, we see that
\begin{align*}(x^2+x+1)^{p-2}\eq&1+\binom{p-2}{1}_2x^{p-1}
+\binom{p-2}{0}_2x^{p-2}
\\&+\sum_{k=2}^{p-2}\left[ \binom{p-2}{k}_2+\binom{p-2}{p-1-k}_2\right] x^{k-1}\pmod p.
\end{align*}
By Lemma 2.1 and \eqref{Tp}, we have
$$3\bi{p-2}1_2\eq2\bi{p-1}0_2\eq2\l(\f p3\r)=2\pmod p.$$
Combining this with \eqref{p-20} and the last paragraph, we obtain
the desired \eqref{x^2}. \qed

\section{Proof of Theorem 1.1}
\setcounter{lemma}{0}
\setcounter{theorem}{0}
\setcounter{corollary}{0}
\setcounter{remark}{0}
\setcounter{equation}{0}

We need the following known lemma \cite[Lemma 10]{K04} on determinants.

\begin{lemma}\label{lemma4}
Let $R$ be a commutative ring with identity, and let $P(x)=\sum_{i=0}^{n-1}a_ix^i\in R[x]$. Then we have
$$\det{\left[P(X_iY_j)\right]}_{1\leqslant i<j\leqslant n}=a_0a_1\cdots a_{n-1} \prod_{1\leqslant i< j\leqslant n}(X_i-X_j)(Y_i-Y_j).$$
\end{lemma}

We also need the following known lemma (cf. \cite[Theorem 1.1]{S-FFA}).

\begin{lemma}\label{lemma5}
Let $p$ be an odd prime. For each $p$-adic integer $x$, let $\{x\}_p$ denote the least nonnegative residue of $x$ modulo $p$. Define
\begin{center}$\mathrm{Inv}_p:=\#\{(i,j):1\leqslant i<j\leqslant p-1 \  and \  \{i^{-1}\} _p>\{j^{-1}\} _p\}, $\end{center}
where $\#S$ denotes the cardinality of a set $S$. Then we have
$$\mathrm{Inv}_p\equiv \frac{p+1}{2}\pmod{2}.$$
\end{lemma}

\medskip
\noindent{\bf Proof of Theorem 1.1}.
Recall that $$D_p(1,1)=\left|(i^2+ij+j^2)^{p-2}\right|_{1\leqslant i,j \leqslant p-1}.$$
By Corollary 2.1, we have
\begin{equation}\label{F}(x^2+x+1)^{p-2}\eq\f23(x^{p-1}-1)+F(x),
\end{equation}
where
$$F(x)=\f53
-\f13x^{p-2}+\sum_{k=2}^{p-2}\l(k\l(\f k3\r)+[3\mid k-1]-\f13\r)x^{k-1}\pmod p.$$
By Fermat's little theorem and \eqref{F}, for any $i,j=1,\ldots,p-1$ we have
$$\f{(i^2+ij+j^2)^{p-2}}{j^{2(p-2)}}=\l(\f{i^2}{j^2}+\f ij+1\r)^{p-2}
\eq F\l(\f ij\r)\pmod p,$$
and hence
$$\l(\f{D_p(1,1)}p\r)=\l(\f{|F(i/j)|_{1\ls i,j\ls p-1}}p\r).$$

By Lemma 3.1, we have
\begin{align*}
\l|F\l(\f ij\r)\r|_{1\leqslant i,j \leqslant p-1}=&\ -\frac{5}{9}\prod_{k=2}^{p-2}\l(k\l(\f k3\r)+[3\mid k-1]-\f13\r) \\
&\ \times\prod_{1\leqslant i< j\leqslant p-1}(i-j)\l(\frac{1}{i}-\frac{1}{j}\r).
\end{align*}
In view of Lemma 3.2,
\begin{equation}\label{i-j}
\begin{aligned}&\prod_{1\leqslant i<j\leqslant p-1}(i-j)\l(\frac{1}{i}-\frac{1}{j}\r)
\\=&\ (-1)^{\mathrm{Inv}_p}\prod_{1\leqslant i<j\leqslant p-1}(i-j)^2=(-1)^{(p+1)/2}\prod_{j=2}^{p-1}((j-1)!)^2.
\end{aligned}
\end{equation}
Thus
\begin{align*}
&\ \left|F\l(\f ij\r)\right|_{1\leqslant i,j \leqslant p-1}\\
\equiv &\ (-1)^{(p+1)/2+1}\cdot \frac{5}{9}\prod_{r=0}^2\prod_{k=2\atop k\eq r\pmod3}^{p-2}\l(k\l(\f k3\r)+[3\mid k-1]-\f13\r)
\times\prod_{j=1}^{p-2}(j!)^2\\
=&\ \f{(-1)^{(p+1)/2}}{(-3)^{(p-1)/3}}\prod_{i=0}^{(p-4)/3}\l((3i+1)+\f23\r)\l(-(3i+2)-\f13\r)
\times\prod_{j=1}^{p-2}(j!)^2
\\
=&\ \f{(-1)^{(p+1)/2+(p-1)/3}}
{3^{p-1}}\prod_{i=0}^{\frac{p-4}{3}}(9i+5)(-9i-7)
\times\prod_{j=1}^{p-2}(j!)^2\pmod{p}.
\end{align*}
As
$$9\l(\frac{p-4}{3}-i\r)+7=3(p-4)+7-9i=3p-5-9i\equiv -(9i+5)\pmod{p}$$
for any $i=0,\ldots,(p-4)/3$, by the above and Lemma 3.2 we have
$$\l|F\l(\f ij\r)\r|_{1\ls i,j\ls p-1}\eq(-1)^{(p+1)/2+(p-1)/3}\prod_{i=0}^{(p-4)/3}(9i+5)^2\times\prod_{j=1}^{p-2}(j!)^2
\pmod p.$$
For $0\ls i\ls(p-4)/3$, clearly $9i+5\ls 3p-7<3p$, and $9i+5\not=p$ since $p\eq1\pmod3$.
Note that
$$9i+5=2p\ \t{for some}\ i=0,\ldots,\f{p-4}3
\iff p\eq-2\eq7\pmod 9.$$
Thus
$$\l(\f{|F(i/j)|_{1\ls i,j\ls p-1}}p\r)
=\l(\f{-1}p\r)^{(p+1)/2+(p-1)/3}\times\begin{cases}0&\t{if}\ p\eq7\pmod 9,
\\1&\t{otherwise}.\end{cases}$$
Therefore
$$\l(\f{D_p(1,1)}p\r)=\l(\f{|F(i/j)|_{1\ls i,j\ls p-1}}p\r)
=\begin{cases}0&\t{if}\ p\eq7\pmod 9,
\\1&\t{otherwise}.\end{cases}$$
This concludes our proof of Theorem 1.1.
\qed

\section{Proof of Theorem 1.2}
\setcounter{lemma}{0}
\setcounter{theorem}{0}
\setcounter{corollary}{0}
\setcounter{remark}{0}
\setcounter{equation}{0}

Let  $A,\ B\in\Z$. The Lucas sequence $u_n=u_n(A,B)\ (n\in\N)$ is defined as follows:
$$u_0=0,\ u_1=1,\ \t{and}\ u_{n+1}=Au_{n}-Bu_{n-1}\ \t{for}\ n=1,2,3,\ldots.$$
Let $\al$ and $\beta$ be the two roots of the quadratic equation $x^2-Ax+B=0$.
By Binet's formula,
$$(\al-\beta)u_n=\al^n-\beta^n\ \quad\t{for all}\ n\in\N.$$

The following lemma is well-known (see, e.g., \cite[Lemma 2.3]{S10}).

\begin{lemma} \label{Lem-AB} Let $A,B\in\Z$, and let $p$ be an odd prime. Then
\begin{equation}\label{up} u_p(A,B)\eq\l(\f{A^2-4B}p\r)\pmod p.
\end{equation}
Provided $p\nmid B$, we also have
 \begin{equation}\label{mod p} u_{p-(\f{A^2-4B}p)}(A,B)\eq0\pmod p.
 \end{equation}
 \end{lemma}

\begin{lemma}\label{Lem-uk} Let $b,c\in\Z$, and let $p$ be an odd prime. Then
\begin{equation}\label{pk-u}\binom{p-1}{p-k}_{b,c}\equiv u_k(-b,c)\pmod{p}
\ \ \t{for all}\ k=0,1,\ldots,p-1.
\end{equation}
\end{lemma}
\Proof. Obviously,
$$\binom{p-1}{p}_{b,c}=0=u_0(-b,c) \ \t{and}\  \binom{p-1}{p-1}_{b,c}=1=u_1(-b,c).$$

Now let $k\in\{1,\ldots,p-2\}$, and assume that $\bi{p-1}{p-j}_{b,c}\eq u_j(b,c)\pmod p$
for all $j=0,\ldots,k$.
By \eqref{recur}, we have
$$
\binom{p}{p-k}_{b,c}=\binom{p-1}{p-k-1}_{b,c}+b\binom{p-1}{p-k}_{b,c}+c\binom{p-1}{p-k+1}_{b,c}.
$$
Since
$$
\binom{p}{p-k}_{b,c}\equiv 0 \pmod{p}
$$
by \eqref{pk}, we have
\begin{align*}\binom{p-1}{p-k-1}_{b,c}&\eq-b\binom{p-1}{p-k}_{b,c}-c\binom{p-1}{p-k+1}_{b,c}
\\&\eq -bu_k(-b,c)-cu_{k-1}(-b,c)=u_{k+1}(-b,c)\pmod p.
\end{align*}

By the above, we have proved \eqref{pk-u} by induction. \qed
\medskip

\begin{lemma}\label{Lem-U}
Let $p$ be an odd prime, and let $b,c\in\Z$ with $p\nmid c(b^2-4c)$. Let
\begin{equation}\label{U}
U(k)= \bi{p-2}k_{b,c}+c^{p-1-k}\bi{p-2}{p-1-k}_{b,c}.
\end{equation}

{\rm (i)} If $U(k)\equiv 0\pmod{p}$ for some $k\in\{2,\ldots,p-2\}$, then
$$
\l(\f{D_p(b,c)}{p}\right)=0.
$$

{\rm (ii)}
If $U(k)\not\equiv 0\pmod{p}$ for all $2\ls k\ls p-2$, then
\begin{equation}\label{D-p}\begin{aligned}
&\left(\f cp\r)^{(p-1)(p-3)/8}\l(\f{D_p(b,c)}{p}\right)\\
=&\left(\frac{4c-b^2+2c(\f{b^2-4c}p)}p\r)
\l(\f{2cu_{p-1}(-b,c)-b}p\r)\l(\f{U(p-2)U(\frac{p-1}{2})}{p}\right).\\
\end{aligned}\end{equation}
\end{lemma}
\Proof.  Recall that
$$D_p(b,c)=\left|(i^2+bij+cj^2)^{p-2}\right|_{1\leqslant i,j \leqslant p-1}.$$
By Lemma 2.2, we have
\begin{equation}\label{G}(x^2+bx+c)^{p-2}\eq\bi{p-2}1_{b,c}(x^{p-1}-1)+G(x),
\end{equation}
where
\begin{align*}
G(x)&=c^{p-2}+\bi{p-2}1_{b,c}+\bi{p-2}0_{b,c}x^{p-2}\\
&\ \ +\sum_{1<k<p-1}\l(\bi{p-2}k_{b,c}+c^{p-1-k}\bi{p-2}{p-1-k}_{b,c}\r)x^{k-1}.\\
\end{align*}

By Fermat's little theorem and \eqref{G}, for any $i,j=1,\ldots,p-1$, we have
$$\f{(i^2+bij+cj^2)^{p-2}}{j^{2(p-2)}}=\l(\f{i^2}{j^2}+\f {bi}j+c\r)^{p-2}
\eq G\l(\f ij\r)\pmod p.$$
Thus
$$\left(\frac{D_p(b,c)}{p}\right)=\left(\frac{|G(i/j)|_{1\leqslant i,j \leqslant p-1}}{p}\right).$$

In view of Lemma \ref{lemma4} and the equality \eqref{i-j},
\begin{align*}
&\l|G\l(\f ij\r)\r|_{1\leqslant i,j\leqslant p-1}\\
=&\left(c^{p-2}+\binom{p-2}{1}_{b,c}\right)\binom{p-2}{0}_{b,c} \prod_{k=2}^{p-2}U(k)\times\prod_{1\le i<j\le p-1}(i-j)\l(\frac{1}{i}-\frac{1}{j}\r)\\
=&\left(c^{p-2}+\binom{p-2}{1}_{b,c}\right)\binom{p-2}{0}_{b,c}U(p-2)U\l(\frac{p-1}{2}\r)\\
&\times\prod_{k=2}^{\frac{p-3}{2}}U(k)U(p-1-k)\times (-1)^{\f {p+1}2}
\prod_{j=1}^{p-2}(j!)^2
\end{align*}
and hence
\begin{equation}\label{DG}\begin{aligned}\l(\f{D_p(b,c)}p\r)=\ &\l(\f{(c^{p-2}+\bi{p-2}1_{b,c})\bi{p-2}0_{b,c}}p\r)
\l(\f{U(p-2)U((p-1)/2)}p\r)
\\&\ \times\prod_{k=2}^{(p-3)/2}\l(\f{U(k)U(p-1-k)}p\r).
\end{aligned}
\end{equation}

(i) If $U(k)\eq0\pmod p$ for some $2\ls k\ls p-2$, then by \eqref{DG} we immediately have
$(\f{D_p(b,c)}p)=0$.

(ii) Now suppose that $U(k)\not\eq0\pmod p$ for any $2\ls k\ls p-2$.
For each $k=1,\ldots,p-1$, with the aid of \eqref{U} we get
\begin{align*}
U(p-k-1)&=c^k\bi{p-2}k_{b,c}+ \bi{p-2}{p-1-k}_{b,c}\\
&\equiv c^k\bi{p-2}k_{b,c}+c^{p-1}\bi{p-2}{p-1-k}_{b,c}\eq c^kU(k)
\pmod p.\end{align*}
Thus
\begin{align*}\prod_{k=2}^{(p-3)/2}\l(\f{U(k)U(p-1-k)}p\r)&=\prod_{k=2}^{(p-3)/2}\l(\f {c^kU(k)^2}p\r)
\\&=\l(\f cp\r)^{\sum_{k=2}^{(p-3)/2}k}=\l(\f cp\r)^{(p-1)(p-3)/8-1},
\end{align*}
and hence \eqref{DG} has the following equivalent form:
\begin{equation}\label{4.7}\begin{aligned}&\l(\f{D_p(b,c)}p\r)\l(\f cp\r)^{(p-1)(p-3)/8}
\\=\ &\l(\f{(1+c\bi{p-2}1_{b,c})\bi{p-2}0_{b,c}}p\r)\l(\f{U(p-2)U((p-1)/2)}p\r).
\end{aligned}
\end{equation}
By Lemmas 2.1 and 4.1, the equality \eqref{-k} and the congruence \eqref{up}, we have
\begin{align*}
&\l(1+c\bi{p-2}1_{b,c}\r)\bi{p-2}0_{b,c}\\
\equiv &\l(1+\f {2c}{4c-b^2}\binom{p-1}{0}_{b,c}\r)\f 1{4c-b^2}\l(\binom{p-1}{-1}_{b,c}+c\binom{p-1}{1}_{b,c}-b\r)\\
= &\l(1+\f {2c}{4c-b^2}\binom{p-1}{0}_{b,c}\r)\f 1{4c-b^2}\l(2c\binom{p-1}{1}_{b,c}-b\r)\\
\equiv &\l(\f 1{4c-b^2}\r)^2\l((4c-b^2)+2cu_{p}(-b,c)\r)\l(2cu_{p-1}(-b,c)-b\r)
\\\eq&\l(\f 1{4c-b^2}\r)^2\l((4c-b^2)+2c\l(\f{b^2-4c}p\r)\r)\l(2cu_{p-1}(-b,c)-b\r) \pmod p.
\end{align*}
Combining this with \eqref{4.7}, we immediately get the desired \eqref{D-p}.

In view of the above, we have completed our proof of Lemma \ref{Lem-U}. \qed

\medskip

\noindent{\bf Proof of Theorem 1.2}.
In view of Binet's formula, for any $k\in\N$
we have $$u_k:=u_k(-2,2)=\f{(-1+i)^k-(-1-i)^k}{2i}$$
and thus
\begin{equation}\label{-22}
u_k=(-4)^{\lfloor \frac{k}{4}\rfloor}
\times\begin{cases}
0       & \t{if}\ k\equiv 0\pmod{4},\\
 1 &\t{if}\ k\equiv 1\pmod{4},  \\
-2   &\t{if}\ k\equiv 2\pmod{4},  \\
 2          &\t{if}\ k\equiv 3\pmod{4}, \end{cases}
  \end{equation}
 which can also be proved easily by induction.
By Lemma \ref{Lem-uk},
\begin{equation}
\label{22} \binom{p-1}{p-k}_{2,2} \equiv u_k\pmod p\ \ \t{for all}\ k=0,1,\ldots,p-1.
\end{equation}

Let $k\in\{2,\ldots,p-2\}$, and
$$U(k)=\bi{p-2}k_{2,2}+2^{p-1-k}\bi{p-2}{p-1-k}_{2,2}.$$
By Lemma \ref{Lem2.1},
$$4\bi{p-2}k_{2,2}\eq(k+1)\bi{p-1}{k-1}_{2,2}-2(k-1)\bi{p-1}{k+1}_{2,2}\pmod p$$
and
\begin{align*}4\bi{p-2}{p-1-k}_{2,2}&\eq(p-k)\bi{p-1}{p-2-k}_{2,2}-2(p-2-k)\bi{p-1}{p-k}_{2,2}
\\&\eq-k\bi{p-1}{p-2-k}_{2,2}+(2k+4)\bi{p-1}{p-k}_{2,2}\pmod p.
\end{align*}
Thus, by the above, we have
\begin{align*}4U(k)&\eq(k+1)\bi{p-1}{k-1}_{2,2}-2(k-1)\bi{p-1}{k+1}_{2,2}
\\&\ \ +2^{p-1-k}\l((2k+4)\bi{p-1}{p-k}_{2,2}-k\bi{p-1}{p-2-k}_{2,2}\r).
\end{align*}
Therefore, with the aid of \eqref{22}, we get
\begin{equation}\label{4U} \begin{aligned}4U(k)&\eq(k+1)u_{p-k+1}-2(k-1)u_{p-k-1}
\\&\ \ +2^{-k}((2k+4)u_k-ku_{k+2})
\pmod p.
\end{aligned}\end{equation}

Now we handle the case $p\eq5\pmod 8$. Applying \eqref{4U} with $k=(p-1)/2$
and noting $(\f 2p)=-1$, we obtain
\begin{align*}
4U\l(\frac{p-1}{2}\r)\eq&\ \frac{p+1}{2}u_{(p+3)/2}
-2\times\frac{p-3}{2}u_{(p+1)/2-1}\\
&\ +2^{-(p-1)/{2}}\left((p+3)u_{(p-1)/2}-\f {p-1}2 u_{(p+3)/2}\right)\\
\equiv&\ \f12u_{(p+3)/{2}}+3u_{(p-1)/{2}}+\l(\f 2p\r)\left(3u_{(p-1)/2}+\f12 u_{(p+3)/2}\right)
\\\eq&\ 0\pmod{p}.
\end{align*}
Combining this with Lemma 4.3,  we see that
$$\left(\frac{D_p(2,2)}{p}\right)=0.$$

Below we assume that $p\eq1\pmod 8$ and write $p=8q+1$ with $q\in\Z^+$.
Let $k\in\{2,\ldots,p-2\}$, and write $k=4s+r$ with $s\in\N$ and $r\in\{0,1,2,3\}$.
We want to show that $U(k)\not\eq0\pmod p$.
\medskip

{\it Case} 1. $r=0$.

In this case, by \eqref{4U} and \eqref{-22} we have
\begin{align*}
4U(k)\equiv &(k+1)u_{p-k+1}-2(k-1)u_{p-k-1}+2^{-k}((2k+4)u_k-ku_{k+2})\\
\equiv& -2(k+1)(-4)^{\lfloor \frac{p-k+1}{4}\rfloor}+2^{-k}\cdot 2k(-4)^{\lfloor \f k4\rfloor}   \\
=& -2(-4)^{2q-s}(k+1)+2^{-4s+1}(-4)^sk   \\
\equiv& -2\l(\f 2p\r)(-4)^{-s}(k+1)+2^{-4s+1}(-4)^sk \\
\equiv& -2(-4)^{-s}(k+1)+2^{-4s+1}(-4)^sk\pmod{p}. \\
\end{align*}
So
$$(-4)^{s+1}U(k)\equiv 2(k+1)-2^{-4s+1}(-4)^{2s}k\eq 2\pmod{p},$$
and hence $U(k)\not\eq0\pmod p$.
\medskip

{\it Case} 2. $r=1$.

In this case, by \eqref{4U} and \eqref{-22} we get
\begin{align*}
4U(k)\equiv&\ (k+1)(-4)^{\lfloor \frac{p-k+1}{4}\rfloor}-4(k-1)(-4)^{\lfloor \frac{p-k-1}{4}\rfloor}\\
&\ +2^{-k}((2k+4)(-4)^{\lfloor \frac{k}{4}\rfloor}-2k(-4)^{\lfloor \frac{k+2}{4}\rfloor}) \\
\equiv&\ (k+1)(-4)^ {2q-s}-4(k-1)(-4)^{2q-s-1} \\
&\ +2^{-4s-1}((2k+4)(-4)^{s}-2k(-4)^{s})   \\
\equiv&\ 2(-4)^{2q-s}k+2^{-4s-1+2}(-4)^{s} \\
\equiv&\ \l(\f 2p\r)2(-4)^{-s}k+2^{-4s+1}(-4)^{s} \\
\equiv&\ 2(-4)^{-s}k+2^{-4s+1}(-4)^{s}\pmod{p},
\end{align*}
and hence
$$
-(-4)^{s+1}U(k)\eq 2k+2^{-4s+1}(-4)^{2s}\equiv 2k+2\not\eq0 \pmod{p}.
$$
Therefore $U(k)\not\eq0\pmod p$.
\medskip

{\it Case} 3. $r=2$.

In light of \eqref{4U} and \eqref{-22}, we have
\begin{align*}
4U(k)\equiv&\ -2(k-1)(-4)^{\lfloor \frac{p-k-1}{4}\rfloor}\times (-2)+2^{-k}(2k+4)(-4)^{\lfloor \frac{k}{4}\rfloor}\times (-2)   \\
\equiv&\ 4(k-1)(-4)^{2q-s-1}-2^{-4s-2}(4k+8)(-4)^{s}   \\
\equiv&\ -\l(\f 2p\r)(-4)^{-s}(k-1)-2^{-4s}(-4)^{s}(k+2) \\
\equiv&\ -(-4)^{-s}(k-1)-2^{-4s}(-4)^{s}(k+2)\pmod{p},
\end{align*}
and hence
$$(-4)^{s+1}U(k)\equiv (k-1)+2^{-4s}(-4)^{2s}(2+k)\equiv 2k+1 \pmod{p}.$$
 Note that $2k+1\not=p$ since $p\eq1\pmod8$ and $k\eq2\pmod 4$. Therefore $U(k)\not\eq0\pmod p$.

\medskip

{\it Case} 4. $r=3$.

By \eqref{4U} and \eqref{-22}, we have
\begin{align*}
4U(k)\equiv&\ 2(k+1)(-4)^{\lfloor \frac{p-k+1}{4}\rfloor}-2(k-1)(-4)^{\lfloor \frac{p-k-1}{4}\rfloor}\\
&\ +2^{-k}(4(2+k)(-4)^{\lfloor \frac{k}{4}\rfloor}-k(-4)^{\lfloor \frac{k+2}{4}\rfloor})   \\
=&\ 2(k+1)(-4)^{2q-s-1}-2(k-1)(-4)^{2q-s-1}\\
&\ +2^{-4s-3}(4(2+k)(-4)^{s}-k(-4)^{s+1})   \\
\equiv&\ -(-4)^{2q-s}+2^{-4s-1}(-4)^{s}(2k+2) \\
\equiv&\ -\l(\f 2p\r)(-4)^{-s}+2^{-4s}(-4)^{s}(k+1) \\
\equiv&\ -(-4)^{-s}+2^{-4s}(-4)^{s}(k+1)\pmod{p}.
\end{align*}
So
$$(-4)^{s+1}U(k)\equiv 1-2^{-4s}(-4)^{2s}(k+1) =-k\pmod{p},$$
and hence $U(k)\not\eq0\pmod p$.

By the above analysis, $U(k)\not\eq0\pmod p$ for each $k=2,3,\dots p-2$.
Note that
$$4\times2^2-2^2+2\times2\l(\f{2^2-4\times2}p\r)=4+4\l(\f{-4}p\r)=8,$$
and $(\f 8p)=(\f 2p)=1$ since $p\eq1\pmod 8$.
Thus, by Lemma \ref{Lem-U}(ii), we have
\begin{equation}\label{dd} \l(\f{D_p(2,2)}p\r)=\l(\f{2u_{p-1}(-2,2)-1}p\r)\l(\f{U(p-2)U((p-1)/2)}p\r).
\end{equation}

Clearly, $$u_{p-1}=u_{p-\l(\f{(-2)^2-4\times2}p\r)}\eq0\pmod p$$ by \eqref{mod p}. Thus
$$\l(\f{2u_{p-1}(-2,2)-1}p\r)=\l(\f{-1}p\r)=1.$$
In view of \eqref{4U} and Lemma \ref{Lem-AB}, we see that
\begin{align*}4U(p-2)&\equiv -u_3+3u_1+2^{p-1}u_p
\\&\eq -2+3\times1+\l(\f{(-2)^2-4\times2}p\r)= 2\pmod{p}
\end{align*}
and
\begin{align*}
4U\l(\frac{p-1}{2}\r)&=\frac{p+1}{2}u_{(p+3)/2}-2\l(\frac{p-1}{2}-1\r)u_{(p-1)/2}\\
&+2^{-(p-1)/2}\left((p+3)u_{(p-1)/2}-\frac{p-1}{2}u_{(p+3)/{2}}\right)\\
&\equiv \f 12 \l(1+\l(\f 2p\r)\r)u_{(p+3)/2}+3\l(1+\l(\f 2p\r)\r)u_{(p-1)/{2}}\\
&\equiv u_{(p+3)/{2}}+6u_{(p-1)/{2}}=-2\times (-4)^{\lfloor \frac{p+1}{8}\rfloor}\pmod{p}.
\end{align*}
(Note that we use \eqref{-22} in the last step.)

Combining the last paragraph with \eqref{dd}, we immediately obtain that
$$
\left(\frac{D_p(2,2)}{p}\right)=1.
$$
This concludes our proof of Theorem 1.2. \qed

\Ack. The authors would like to thank the anonymous referee for helpful comments.
This work was supported by the National Natural Science Foundation of China (grant no. 11971222).

\setcounter{conjecture}{0}
\end{document}